# ON THE TWO-TIMES DIFFERENTIABILITY OF THE VALUE FUNCTIONS IN THE PROBLEM OF OPTIMAL INVESTMENT IN INCOMPLETE MARKETS

By Dmitry Kramkov[1] and Mihai Sîrbu

*Carnegie Mellon University and Columbia University*

We study the two-times differentiability of the value functions of the primal and dual optimization problems that appear in the setting of expected utility maximization in incomplete markets. We also study the differentiability of the solutions to these problems with respect to their initial values. We show that the key conditions for the results to hold true are that the relative risk aversion coefficient of the utility function is uniformly bounded away from zero and infinity, and that the prices of traded securities are sigma-bounded under the numéraire given by the optimal wealth process.

**1. Introduction and main results.** We study a similar financial framework to the one in [7] and refer to this paper for more details and references. We consider a model of a security market which consists of $d + 1$ assets, one bond and $d$ stocks. We work in discounted terms, that is, we suppose that the price of the bond is constant and denote by $S = (S^i)_{1 \leq i \leq d}$ the price process of the $d$ stocks. The process $S$ is assumed to be a semimartingale on a filtered probability space $(\Omega, \mathcal{F}, (\mathcal{F}_t)_{0 \leq t \leq T}, \mathbb{P})$. Here $T$ is a finite time horizon and $\mathcal{F} = \mathcal{F}_T$.

A (self-financing) portfolio is defined as a pair $(x, H)$, where the constant $x$ represents the initial capital and $H = (H^i)_{1 \leq i \leq d}$ is a predictable $S$-integrable process, $H^i_t$ specifying how many units of asset $i$ are held in the portfolio at time $t$. The wealth process $X = (X_t)_{0 \leq t \leq T}$ of the portfolio evolves in time as the stochastic integral of $H$ with respect to $S$:

$$(1) \qquad X_t = X_0 + \int_0^t H_u \, dS_u, \qquad 0 \leq t \leq T.$$

Received February 2005; revised December 2005.

[1]Supported in part by the NSF Grants DMS-01-39911 and DMS-05-05414.

We denote by $\mathcal{X}(x)$ the family of nonnegative wealth processes with initial value $x$:

(2) $$\mathcal{X}(x) = \{X \geq 0 \colon X \text{ is defined by (1) with } X_0 = x\}.$$

A probability measure $\mathbb{Q} \sim \mathbb{P}$ is called an *equivalent local martingale measure* if any $X \in \mathcal{X}(1)$ is a local martingale under $\mathbb{Q}$. The family of equivalent local martingale measures is denoted by $\mathcal{Q}$. We assume throughout that

(3) $$\mathcal{Q} \neq \varnothing.$$

This condition is intimately related to the absence of arbitrage opportunities on the security market. See [1] for precise statements and references.

We also consider an economic agent in our model, whose preferences over terminal wealth are modeled by a utility function $U = (U(x))_{x>0}$. The function $U$ is assumed to be strictly concave, strictly increasing and continuously differentiable and to satisfy the Inada conditions:

(4) $$U'(0) = \lim_{x \to 0} U'(x) = \infty, \qquad U'(\infty) = \lim_{x \to \infty} U'(x) = 0.$$

In what follows we set $U(0) = \lim_{x \to 0} U(x)$ and $U(x) = -\infty$ for all $x < 0$.

For a given initial capital $x > 0$, the goal of the agent is *to maximize the expected value of terminal utility*. The value function of this problem is denoted by

(5) $$u(x) = \sup_{X \in \mathcal{X}(x)} \mathbb{E}[U(X_T)].$$

Intuitively speaking, the value function $u = (u(x))_{x>0}$ plays the role of the utility function of the investor at time 0, if he/she subsequently invests in an optimal way. To exclude the trivial case, we assume that $u$ is finite:

(6) $$u(x) < \infty, \qquad x > 0.$$

A well-known tool in studying the optimization problem (5) is the use of duality relationships in the spaces of convex functions and semimartingales. Following [7], we define the dual optimization problem to (5) as follows:

(7) $$v(y) = \inf_{Y \in \mathcal{Y}(y)} \mathbb{E}[V(Y_T)], \qquad y > 0.$$

Here $V$ is the convex conjugate function to $U$, that is,

$$V(y) = \sup_{x>0}\{U(x) - xy\}, \qquad y > 0,$$

and $\mathcal{Y}(y)$ is the family of nonnegative supermartingales $Y$ that are dual to $\mathcal{X}(1)$ in the following sense:

(8) $$\mathcal{Y}(y) = \{Y \geq 0 \colon Y_0 = y \text{ and } XY \text{ is a supermartingale for all } X \in \mathcal{X}(1)\}.$$

Note that the set $\mathcal{Y}(1)$ contains the density processes of all $\mathbb{Q} \in \mathcal{Q}$.

The optimization problems (5) and (7) are well studied. For example, it was shown in [7] that the value functions $u$ and $v$ are conjugate, that is,

$$(9) \qquad v(y) = \sup_{x>0}\{u(x) - xy\}, \qquad y > 0,$$

and that the *minimal market independent condition* on $U$ that implies the continuous differentiability of the value functions $u$ and $v$ on $(0, \infty)$ and the existence of the solutions $X(x)$ and $Y(y)$ to (5) and (7) for all $x > 0$ and $y > 0$ is that the asymptotic elasticity of $U$ is strictly less than 1, that is,

$$(10) \qquad \limsup_{x\to\infty} \frac{xU'(x)}{U(x)} < 1.$$

In addition, if $y = u'(x)$ then

$$(11) \qquad U'(X_T(x)) = Y_T(y)$$

and the product $X(x)Y(y)$ is a martingale. Hereafter, we shall use these results without further comment.

In this paper we are interested in the existence and the computation of the *second derivatives* $u''(x)$ and $v''(y)$ of the value functions $u$ and $v$ and of the *first derivatives* $X'(x)$ and $Y'(y)$ of the solutions $X(x)$ and $Y(y)$ with respect to $x$ and $y$. In addition to a purely theoretical interest (see Remark 1 below) our study of these questions is also motivated by some applications. For example, in [9] we use the results of the current paper to perform the sensitivity analysis of utility-based prices with respect to the number of nontraded contingent claims.

REMARK 1. If $S$ is a Markov diffusion process, then (5) becomes a typical stochastic control problem and can be studied using PDE methods. In this case, the two-times differentiability of the value function $u$ is closely related to the existence of the *classical* solution of the corresponding Bellman equation. We refer to [10] for the deep treatment of this topic.

To give positive answers to the above questions we need to impose additional conditions on the utility function $U$ and the price process $S$. The conditions on $U$ are stated in the following assumption:

ASSUMPTION 1. The utility function $U$ is two-times differentiable on $(0, \infty)$ and its relative risk aversion coefficient

$$(12) \qquad A(x) = -\frac{xU''(x)}{U'(x)}, \qquad x > 0,$$

is uniformly bounded away from zero and infinity, that is, there are constants $c_1 > 0$ and $c_2 < \infty$ such that

$$(13) \qquad c_1 < A(x) < c_2, \qquad x > 0.$$

In Lemma **3** below we prove that the bounds **(13)** on the relative risk aversion coefficient imply both the Inada conditions **(4)** and the condition **(10)** on the asymptotic elasticity. Note that $U$ is two-times differentiable at $x > 0$ and $U''(x) < 0$ if and only if the conjugate function $V$ is two-times differentiable at $y = U'(x)$, and that in this case

$$(14) \qquad V''(y) = -\frac{1}{U''(x)}.$$

It follows that **(13)** is equivalent to the following condition:

$$(15) \qquad \frac{1}{c_2} < B(y) < \frac{1}{c_1}, \qquad y > 0,$$

where

$$(16) \qquad B(y) = -\frac{yV''(y)}{V'(y)}, \qquad y > 0.$$

Note that if $y = U'(x)$ then

$$(17) \qquad B(y) = \frac{1}{A(x)}$$

is the *relative risk tolerance coefficient* of $U$ computed at $x$.

To facilitate the formulation of the assumption on the price process $S$ we give the following definition:

DEFINITION 1. A $d$-dimensional semimartingale $R$ is called *sigma-bounded* if there is a strictly positive predictable (one-dimensional) process $h$ such that the stochastic integral $\int h\,dR$ is well defined and is locally bounded.

This definition has been motivated by a similar concept of sigma-martingales which plays the key role in the fundamental theorem of asset pricing for unbounded processes; see [1] for the details. Both these notions are instances of a more general concept of sigma-localization studied in [6]. Appendix **A** below contains a number of equivalent reformulations of Definition **1**, as well as other useful results on sigma-bounded semimartingales and martingales.

In principle, any strictly positive wealth process $X$ can be chosen as a new currency or a *numéraire* in the model. In this case, the $(d+1)$-dimensional semimartingale

$$(18) \qquad S^X = \left(\frac{1}{X}, \frac{S}{X}\right)$$

has the economic interpretation as the prices of traded securities (the bond and the stocks) discounted by $X$. It is well known that a process $\widetilde{X}$ is a stochastic integral with respect to $S^X$ if and only if $\widetilde{X}X$ is a stochastic integral with respect to $S$. In other words, $\widetilde{X}$ is a wealth process under the numéraire $X$ if and only if $\widetilde{X}X$ is a wealth process in the original financial model, where the role of the money denominator is played by the bond. We shall use this fact on several occasions.

Hereafter, we fix an initial capital $x > 0$ and denote $y = u'(x)$. The following assumption plays a key role in the proofs of our main results. It states that the prices of traded securities are sigma-bounded under the numéraire given as the solution $X(x)$ to (5).

Assumption 2. The price process of the traded securities discounted by $X(x)$, that is, the $(d+1)$-dimensional semimartingale

$$(19) \qquad S^{X(x)} = \left( \frac{1}{X(x)}, \frac{S}{X(x)} \right)$$

is sigma-bounded.

The direct verification of Assumption 2 is feasible only if we can compute the optimal wealth process $X(x)$ explicitly. In other cases, some sufficient "qualitative" conditions on the financial model should be used. An example of such a condition is given by Theorem 3 in Appendix A. This theorem shows that *all* semimartingales defined on the filtered probability space $(\Omega, \mathcal{F}, (\mathcal{F}_t)_{0 \le t \le T}, \mathbb{P})$ are sigma-bounded (and, hence, Assumption 2 holds true trivially) if the family of purely discontinuous martingales admits a finite-dimensional basis (from the point of view of stochastic integration). In particular, as Proposition 2 in Appendix A shows, all semimartingales are sigma-bounded if the financial model is *complete* (or can be extended to become complete by adding to it a finite number of additional securities). It is interesting to note that there are complete financial models where the "locally bounded" version of Assumption 2 fails to hold true; see Example 5 in Section 4.

To state our main results we also need to define two auxiliary optimization problems. Let $\mathbb{R}(x)$ be the probability measure whose Radon–Nikodym derivative under $\mathbb{P}$ is given by

$$(20) \qquad \frac{d\mathbb{R}(x)}{d\mathbb{P}} = \frac{X_T(x)Y_T(y)}{xy} = \frac{X_T(x)U'(X_T(x))}{xu'(x)}.$$

(A similar change of measure was used in the particular case $U(x) = -x^2$ in [2].) Let $\mathbf{H}_0^2(\mathbb{R}(x))$ be the space of square integrable martingales under $\mathbb{R}(x)$ with initial value 0. We would like to point out that square integrability is naturally related to the existence of second-order derivatives, as the reader can see below. Denote by $\mathcal{M}^2(x)$ the subspace of $\mathbf{H}_0^2(\mathbb{R}(x))$ that consists of stochastic integrals with respect to $S^{X(x)}$, that is,

$$\mathcal{M}^2(x) = \left\{ M \in \mathbf{H}_0^2(\mathbb{R}(x)) : M = \int H\, dS^{X(x)} \text{ for some } H \right\}.$$

Further, let $\mathcal{N}^2(y) = \mathcal{N}^2(u'(x))$ be the orthogonal complement of $\mathcal{M}^2(x)$ in $\mathbf{H}_0^2(\mathbb{R}(x))$. In other words, $N \in \mathcal{N}^2(y)$ if and only if $N \in \mathbf{H}_0^2(\mathbb{R}(x))$ and $MN$ is a martingale under $\mathbb{R}(x)$ for all $M \in \mathcal{M}^2(x)$.

After these preparations, we formulate the following optimization problems:

$$(21) \qquad a(x) = \inf_{M \in \mathcal{M}^2(x)} \mathbb{E}_{\mathbb{R}(x)}[A(X_T(x))(1 + M_T)^2],$$

$$(22) \qquad b(y) = \inf_{N \in \mathcal{N}^2(y)} \mathbb{E}_{\mathbb{R}(x)}[B(Y_T(y))(1 + N_T)^2], \qquad y = u'(x),$$

where the functions $A$ and $B$ are defined in (12) and (16), respectively. The basic properties of these optimization problems are stated in Lemma 1 below. The proof of this lemma will be given in Section 3 and will follow from its "abstract version," Lemma 2.

LEMMA 1. *Assume that conditions* (3), (6) *and Assumption* 1 *hold true. Let* $x > 0$ *and* $y = u'(x)$. *Then:*

1. *The value functions* $a(x)$ *and* $b(y)$ *defined in* (21) *and* (22) *satisfy*

$$a(x)b(y) = 1$$

   *and*

$$c_1 < a(x) < c_2, \qquad \frac{1}{c_2} < b(y) < \frac{1}{c_1},$$

   *where the constants* $c_1$ *and* $c_2$ *appear in* (13).

2. *The solutions* $M(x)$ *to* (21) *and* $N(y)$ *to* (22) *exist and are unique. In addition,*

$$(23) \qquad A(X_T(x))(1 + M_T(x)) = a(x)(1 + N_T(y)).$$

For $x > 0$ and $y = u'(x)$ we define the semimartingales $X'(x)$ and $Y'(y)$:

$$(24) \qquad X'(x) = \frac{X(x)}{x}(1 + M(x)),$$

$$(25) \qquad Y'(y) = \frac{Y(y)}{y}(1 + N(y)),$$

where $M(x)$ and $N(y)$ are the solutions to (21) and (22), respectively. Note that as $M(x)$ is a stochastic integral with respect to $S^{X(x)}$, the semimartingale $X'(x)$ is a stochastic integral with respect to $S$. In other words, $X'(x)$ is a wealth process.

The following theorem is the main result of our paper:

THEOREM 1. *Let* $x > 0$ *and denote* $y = u'(x)$. *Assume that conditions* (3) *and* (6) *and Assumptions* 1 *and* 2 *hold true. Then:*

1. *The second derivatives $u''(x)$ and $v''(y)$ of the value functions $u$ and $v$ defined in* (5) *and* (7) *exist at $x$ and $y$, respectively, and*

$$(26) \qquad c_1 < -\frac{xu''(x)}{u'(x)} = a(x) < c_2,$$

$$(27) \qquad \frac{1}{c_2} < -\frac{yv''(y)}{v'(y)} = b(y) < \frac{1}{c_1},$$

   *where $a(x)$ and $b(y)$ are the value functions defined in* (21) *and* (22), *and the constants $c_1$ and $c_2$ are given in* (13).

2. *The derivatives of the terminal values of the solutions $X(x)$ to* (5) *and $Y(y)$ to* (7) *with respect to $x$ and $y$ exist and equal the terminal values of the semimartingales $X'(x)$ and $Y'(y)$ defined in* (24) *and* (25). *That is,*

$$(28) \qquad \lim_{\varepsilon \to 0}\left(\frac{X_T(x+\varepsilon) - X_T(x)}{\varepsilon}\right) = X_T'(x),$$

$$(29) \qquad \lim_{\varepsilon \to 0}\left(\frac{Y_T(y+\varepsilon) - Y_T(y)}{\varepsilon}\right) = Y_T'(y),$$

   *where the convergence takes place in $\mathbb{P}$-probability.*

3. *We have*

$$(30) \qquad u''(x) = \mathbb{E}[U''(X_T(x))(X_T'(x))^2],$$

$$(31) \qquad v''(y) = \mathbb{E}[V''(Y_T(y))(Y_T'(y))^2].$$

   *Moreover,*

$$(32) \qquad U''(X_T(x))X_T'(x) = u''(x)Y_T'(y),$$

   *and the products $X(x)Y'(y)$, $X'(x)Y(y)$ and $X'(x)Y'(y)$ are martingales under $\mathbb{P}$.*

REMARK 2. According to (28), the wealth process $X'(x)$ describes how the agent invests a small additional unit of capital. Note that, while the construction of the processes $X'(x)$ and $Y'(y)$ by (24) and (25) was based on Lemma 1 and, hence, did not require Assumption 2, the equalities (28) and (29) hold true, in general, only if this sigma-boundedness assumption is satisfied. Note also that (under the conditions of Theorem 1) equalities (28) and (29) and the martingale properties of $X'(x)Y(y)$ and $X(x)Y'(y)$, where $y = u'(x)$, allow us to compute $X'(x)$ and $Y'(y)$ directly from $X(x)$ and $Y(y)$ without relying on the optimization problems (21) and (22). Finally, we point out that, contrary to a naive conjecture, the processes $X'(x)$ and $Y'(y)$ might not be positive; see Example 4 in Section 4 below.

The proof of Theorem 1 will be given in Section 3 after we study the "abstract" version of this theorem in Section 2. In Section 4 we construct counterexamples to some natural but false conjectures related to our main results. In particular, these counterexamples show that, in general, the upper and the lower bounds in (13), as well as Assumption 2, cannot be removed without affecting the existence of the second derivatives $u''$ and $v''$. Finally, in Appendix A we present some results on sigma-bounded semimartingales and martingales. In particular, Theorem 3 there contains convenient sufficient conditions on the underlying filtered probability space $(\Omega, \mathcal{F}, (\mathcal{F}_t)_{0 \le t \le T}, \mathbb{P})$ that ensure *each* semimartingale $X$ is sigma-bounded and, hence, guarantee the validity of Assumption 2.

**2. The abstract versions of the main results.**   Hereafter, we use the standard notation $\mathbf{L}^0$ for the set of all random variables and $\mathbf{L}^\infty$ for the set of bounded random variables on $(\Omega, \mathcal{F}, \mathbb{P})$. If $\mathbb{Q} \sim \mathbb{P}$, then we denote

$$(33) \qquad \mathbf{L}_0^2(\mathbb{Q}) = \{g \in \mathbf{L}^0 : \mathbb{E}_{\mathbb{Q}}[g] = 0 \text{ and } \mathbb{E}_{\mathbb{Q}}[g^2] < \infty\}.$$

We start with the abstract version of Lemma 1. Let $\mathcal{A}$ and $\mathcal{B}$ be nonempty complementary linear subspaces of $\mathbf{L}_0^2(\mathbb{P})$, that is,

$$(34) \qquad \begin{aligned} \alpha \in \mathcal{A} &\iff \alpha \in \mathbf{L}_0^2(\mathbb{P}) \quad \text{and} \quad \mathbb{E}[\alpha\beta] = 0 \qquad \forall \beta \in \mathcal{B}, \\ \beta \in \mathcal{B} &\iff \beta \in \mathbf{L}_0^2(\mathbb{P}) \quad \text{and} \quad \mathbb{E}[\alpha\beta] = 0 \qquad \forall \alpha \in \mathcal{A}. \end{aligned}$$

Let $\zeta$ be a random variable such that

$$(35) \qquad c_1 < \zeta < c_2,$$

for some constants $0 < c_1 < c_2 < \infty$, and let $\eta$ be the reciprocal to $\zeta$:

$$\eta = \frac{1}{\zeta}.$$

We consider the following optimization problems:

$$(36) \qquad a = \inf_{\alpha \in \mathcal{A}} \mathbb{E}[\zeta(1+\alpha)^2],$$

$$(37) \qquad b = \inf_{\beta \in \mathcal{B}} \mathbb{E}[\eta(1+\beta)^2].$$

LEMMA 2.   *Assume* (34) *and* (35). *Then:*

1. *The numbers $a$ and $b$ defined in* (36) *and* (37) *satisfy*

$$(38) \qquad ab = 1$$

   *and*

$$(39) \qquad c_1 < a < c_2, \qquad \frac{1}{c_2} < b < \frac{1}{c_1},$$

   *where the constants $c_1$ and $c_2$ appear in* (35).

2. *The solutions $\widehat{\alpha}$ to (36) and $\widehat{\beta}$ to (37) exist and are unique. In addition,*

$$\zeta(1+\widehat{\alpha}) = a(1+\widehat{\beta}). \tag{40}$$

PROOF. First, inequalities (39) are easy consequences of (35). Further, let $(\gamma_n)_{n\geq 1}$ be a sequence in $\mathcal{A}$ such that

$$\lim_{n\to\infty} \mathbb{E}[\zeta(1+\gamma_n)^2] = a.$$

As $\zeta \geq c_1$, the sequence $(\gamma_n)_{n\geq 1}$ is bounded in $\mathbf{L}_0^2(\mathbb{P})$. Hence, there is a sequence of convex combinations

$$\alpha_n \in \mathrm{conv}(\gamma_n, \gamma_{n+1}, \dots), \qquad n \geq 1,$$

that converges in $\mathbf{L}_0^2(\mathbb{P})$ to some $\widehat{\alpha}$. As $\mathcal{A}$ is closed in $\mathbf{L}_0^2(\mathbb{P})$, we have $\widehat{\alpha} \in \mathcal{A}$. The convexity of $(1+x)^2$ and the inequality $\zeta \leq c_2$ now imply that

$$\mathbb{E}[\zeta(1+\widehat{\alpha})^2] = \lim_{n\to\infty} \mathbb{E}[\zeta(1+\alpha_n)^2] \leq \limsup_{n\to\infty} \mathbb{E}[\zeta(1+\gamma_n)^2] = a.$$

This proves that $\widehat{\alpha}$ is a solution to (36). The fact that $\widehat{\alpha}$ is the only solution to (36) follows from the strict convexity of $(1+x)^2$.

Using standard arguments from the calculus of variations, we deduce that the optimality of $\widehat{\alpha}$ implies that for any $\alpha \in \mathcal{A}$,

$$\mathbb{E}[\zeta(1+\widehat{\alpha})\alpha] = 0.$$

From the complementary relations (34) between $\mathcal{A}$ and $\mathcal{B}$, we deduce the existence of a constant $c$ and a random variable $\gamma \in \mathcal{B}$ such that

$$\zeta(1+\widehat{\alpha}) = c + \gamma.$$

If we multiply both sides of this equality by $1+\widehat{\alpha}$ and compute the expected value, we get $c = a \neq 0$. Hence, denoting $\beta = \gamma/a$, we deduce that

$$\zeta(1+\widehat{\alpha}) = a(1+\beta), \tag{41}$$

where $\beta \in \mathcal{B}$.

Repeating the same arguments for the optimization problem (37), we deduce the existence and the uniqueness of the solution $\widehat{\beta}$ to this problem, as well as the representation

$$\eta(1+\widehat{\beta}) = b(1+\alpha), \tag{42}$$

for some $\alpha \in \mathcal{A}$.

By multiplying (41) and (42), and using the fact that $\zeta\eta = 1$, we arrive at the equality

$$(1+\widehat{\alpha})(1+\widehat{\beta}) = ab(1+\alpha)(1+\beta),$$

which, after taking the expectation under $\mathbb{P}$, implies the relation (38) between $a$ and $b$. The equality of $\widehat{\beta}$ defined in (41) to $\widehat{\beta}$ now follows from the uniqueness of the solution to (37) and the computations below:

$$\mathbb{E}[\eta(1+\beta)^2] = \mathbb{E}\left[\eta\frac{\zeta^2}{a^2}(1+\widehat{\alpha})^2\right] = \frac{1}{a^2}\mathbb{E}[\zeta(1+\widehat{\alpha})^2] = \frac{1}{a} = b. \qquad \square$$

In the remaining part of this section we state and prove the abstract version of Theorem 1. Let $\mathcal{C}$ and $\mathcal{D}$ be nonempty sets of nonnegative random variables such that:

1. The set $\mathcal{C}$ is bounded in $\mathbf{L}^0$ and contains the constant function $g = 1$:

(43) $$\lim_{n\to\infty}\sup_{g\in\mathcal{C}}\mathbb{P}[|g|\geq n] = 0,$$

(44) $$1 \in \mathcal{C}.$$

2. The sets $\mathcal{C}$ and $\mathcal{D}$ satisfy the bipolar relations

(45) $$\begin{aligned} g \in \mathcal{C} &\iff \mathbb{E}[gh] \leq 1 \qquad \forall\, h \in \mathcal{D}, \\ h \in \mathcal{D} &\iff \mathbb{E}[gh] \leq 1 \qquad \forall\, g \in \mathcal{C}. \end{aligned}$$

For $x > 0$ and $y > 0$, we define the sets

$$\mathcal{C}(x) = x\mathcal{C} = \{xg : g \in \mathcal{C}\},$$
$$\mathcal{D}(y) = y\mathcal{D} = \{yh : h \in \mathcal{D}\},$$

and the optimization problems

(46) $$u(x) = \sup_{g\in\mathcal{C}(x)}\mathbb{E}[U(g)],$$

(47) $$v(y) = \inf_{h\in\mathcal{D}(y)}\mathbb{E}[V(h)].$$

Here $U$ and $V$ are the functions defined in Section 1.

Hereafter, we assume that

(48) $$u(x) < \infty, \qquad x > 0,$$

and that Assumption 1 holds true. From Lemma 3 below and Theorem 3.2 in [7] we deduce that the value functions $u$ and $v$ defined in (46) and (47) are conjugate, that is, (9) holds true, $u$ and $v$ are continuously differentiable on $(0, +\infty)$ and the solutions $g(x)$ to (46) and $h(y)$ to (47) exist for all $x > 0$ and $y > 0$. In addition, if $y = u'(x)$, then

(49) $$U'(g(x)) = h(y),$$

(50) $$\mathbb{E}[g(x)h(y)] = xy.$$

Hereafter, we shall use these results without further comment. As in the previous section, we are interested in the existence of second derivatives $u''(x)$ and $v''(y)$ of the value functions and the first derivatives $g'(x)$ and $h'(y)$ of the solutions to these problems.

For $x > 0$ we denote by $\mathcal{A}^\infty(x)$ the family of bounded random variables $\alpha$ such that $g(x)(1 + c\alpha)$ and $g(x)(1 - c\alpha)$ belong to $\mathcal{C}(x)$ for some $c = c(\alpha) > 0$, where $g(x)$ is the solution to (46). In other words,

$$(51) \qquad \mathcal{A}^\infty(x) = \{\alpha \in \mathbf{L}^\infty : g(x)(1 \pm c\alpha) \in \mathcal{C}(x) \text{ for some } c > 0\}.$$

Similarly, for $y > 0$ we denote

$$(52) \qquad \mathcal{B}^\infty(y) = \{\beta \in \mathbf{L}^\infty : h(y)(1 \pm c\beta) \in \mathcal{D}(y) \text{ for some } c > 0\},$$

where $h(y)$ is the solution to (47).

Hereafter, we fix $x > 0$ and denote $y = u'(x)$. Let $\mathbb{R}(x)$ be the probability measure on $(\Omega, \mathcal{F})$ whose Radon–Nikodym derivative under $\mathbb{P}$ is given by

$$(53) \qquad\qquad \frac{d\mathbb{R}(x)}{d\mathbb{P}} = \frac{g(x)h(y)}{xy}.$$

From the bipolar relations (45) for the sets $\mathcal{C}$ and $\mathcal{D}$ we deduce that the sets $\mathcal{A}^\infty(x)$ and $\mathcal{B}^\infty(y) = \mathcal{B}^\infty(u'(x))$ are orthogonal linear subspaces in $\mathbf{L}_0^2(\mathbb{R}(x))$. That is,

$$(54) \qquad \mathbb{E}_{\mathbb{R}(x)}[\alpha] = \mathbb{E}_{\mathbb{R}(x)}[\beta] = \mathbb{E}_{\mathbb{R}(x)}[\alpha\beta] = 0 \qquad \forall \alpha \in \mathcal{A}^\infty(x), \beta \in \mathcal{B}^\infty(y).$$

Denote by $\mathcal{A}^2(x)$ and $\mathcal{B}^2(y)$ the respective closures of $\mathcal{A}^\infty(x)$ and $\mathcal{B}^\infty(y)$ in $\mathbf{L}_0^2(\mathbb{R}(x))$. From (54) we deduce that $\mathcal{A}^2(x)$ and $\mathcal{B}^2(y)$ are closed *orthogonal* linear subspaces in $\mathbf{L}_0^2(\mathbb{R}(x))$. It turns out that the two-times differentiability of $u$ at $x$ and of $v$ at $y$ depends crucially on the condition that these two subspaces are *complementary* to each other.

ASSUMPTION 3.   The sets $\mathcal{A}^2(x)$ and $\mathcal{B}^2(y)$, where $y = u'(x)$, are complementary linear subspaces in $\mathbf{L}_0^2(\mathbb{R}(x))$. That is,

$$
\begin{aligned}
(55) \qquad & \alpha \in \mathcal{A}^2(x) \quad \iff \quad \alpha \in \mathbf{L}_0^2(\mathbb{R}(x)) \quad \text{and} \quad \mathbb{E}_{\mathbb{R}(x)}[\alpha\beta] = 0 \\
& \hspace{10.5cm} \forall \beta \in \mathcal{B}^2(y), \\
& \beta \in \mathcal{B}^2(y) \quad \iff \quad \beta \in \mathbf{L}_0^2(\mathbb{R}(x)) \quad \text{and} \quad \mathbb{E}_{\mathbb{R}(x)}[\alpha\beta] = 0 \\
& \hspace{10.5cm} \forall \alpha \in \mathcal{A}^2(x).
\end{aligned}
$$

As we show in Section 3, Assumption 3 is the "abstract version" of Assumption 2.

Consider now the optimization problems

$$(56) \qquad a(x) = \inf_{\alpha \in \mathcal{A}^2(x)} \mathbb{E}_{\mathbb{R}(x)}[A(g(x))(1+\alpha)^2],$$

$$(57) \qquad b(y) = \inf_{\beta \in \mathcal{B}^2(y)} \mathbb{E}_{\mathbb{R}(x)}[B(h(y))(1+\beta)^2], \qquad y = u'(x),$$

where the functions $A$ and $B$ are defined in (12) and (16). From Lemma 2 we deduce that if Assumptions 1 and 3 hold true, then the solutions $\alpha(x)$ to (56) and $\beta(y)$ to (57) exist and are unique, and [recalling that $y = u'(x)$]

$$(58) \qquad a(x)b(y) = 1,$$

$$(59) \qquad A(g(x))(1+\alpha(x)) = a(x)(1+\beta(y)).$$

Using this notation, we define the random variables

$$(60) \qquad g'(x) = \frac{g(x)}{x}(1+\alpha(x)),$$

$$(61) \qquad h'(y) = \frac{h(y)}{y}(1+\beta(y)).$$

THEOREM 2. *Let* $x > 0$ *and denote* $y = u'(x)$. *Assume that conditions* (43)–(45), (48) *and Assumptions* 1 *and* 3 *hold true. Then:*

1. *The value functions* $u$ *and* $v$ *defined in* (46) *and* (47) *are two-times differentiable at* $x$ *and* $y$, *respectively, and*

$$(62) \qquad c_1 < -\frac{xu''(x)}{u'(x)} = a(x) < c_2,$$

$$(63) \qquad \frac{1}{c_2} < -\frac{yv''(y)}{v'(y)} = b(y) < \frac{1}{c_1},$$

*where* $a(x)$ *and* $b(y)$ *are defined in* (56) *and* (57) *and the constants* $c_1$ *and* $c_2$ *are given in* (13).

2. *The derivatives of the solutions* $g(x)$ *to* (46) *and* $h(y)$ *to* (47) *with respect to* $x$ *and* $y$ *exist and equal* $g'(x)$ *and* $h'(y)$, *respectively, as defined in* (60) *and* (61), *that is,*

$$(64) \qquad \lim_{\varepsilon \to 0}\left(\frac{g(x+\varepsilon)-g(x)}{\varepsilon}\right) = g'(x),$$

$$(65) \qquad \lim_{\varepsilon \to 0}\left(\frac{h(y+\varepsilon)-h(y)}{\varepsilon}\right) = h'(y),$$

*where the convergence takes place in* $\mathbb{P}$-*probability.*

3. *We have*

$$(66) \qquad u''(x) = \mathbb{E}[U''(g(x))(g'(x))^2],$$

$$(67) \qquad v''(y) = \mathbb{E}[V''(h(y))(h'(y))^2].$$

*Moreover,*

$$(68) \qquad U''(g(x))g'(x) = u''(x)h'(y).$$

Before proceeding to the proof of the theorem, we state two technical results related to the condition (13) on the utility function $U$.

LEMMA 3.   *Assume that the utility function $U = U(x)$ satisfies* (13). *Then the following assertions hold true:*

1. *For any constant $a > 1$ there are constants $0 < b_1 < b_2 < 1$ such that*

$$(69) \qquad b_1 U'(x) < U'(ax) < b_2 U'(x), \qquad x > 0.$$

2. *The function $U$ satisfies the Inada conditions* (4).
3. *The asymptotic elasticity of $U$ is strictly less than 1, that is,* (10) *holds true.*

PROOF.   Let $c_1$ and $c_2$ be the constants defined in (13). Without any loss of generality we can assume that $a > 1$ is sufficiently close to 1 so that

$$1 - c_2 \ln a > 0.$$

Using (13) we deduce

$$U'(x) - U'(ax) = \int_1^a -xU''(tx)\,dt \geq c_1 \int_1^a \frac{U'(tx)}{t}\,dt > c_1 U'(ax)\ln a,$$

$$U'(x) - U'(ax) = \int_1^a -xU''(tx)\,dt \leq c_2 \int_1^a \frac{U'(tx)}{t}\,dt < c_2 U'(x)\ln a.$$

The inequalities in (69) now follow, with

$$b_1 = 1 - c_2 \ln a, \qquad b_2 = \frac{1}{1 + c_1 \ln a}.$$

The Inada conditions (4) follow from (69):

$$U'(\infty) = \lim_{n \to \infty} U'(a^n) \leq \lim_{n \to \infty} (b_2)^n U'(1) = 0,$$

$$U'(0) = \lim_{n \to \infty} U'\left(\frac{1}{a^n}\right) \geq \lim_{n \to \infty} \left(\frac{1}{b_2}\right)^n U'(1) = \infty.$$

Finally, the fact that the second inequality in (69) implies (10) has been proved in [7], Lemma 6.5.   □

LEMMA 4. *Assume that the utility function $U = U(x)$ satisfies* (13). *Let $\zeta$ be a strictly positive random variable such that*

$$\mathbb{E}[|U(\zeta)|] < \infty \tag{70}$$

*and $\eta$ a random variable such that*

$$|\eta| \le K\zeta$$

*for some $K > 0$. Then the function*

$$w(s) = \mathbb{E}[U(\zeta + s\eta)]$$

*is well defined and two-times differentiable for $|s| < \frac{1}{K}$. Furthermore,*

$$w'(s) = \mathbb{E}[U'(\zeta + s\eta)\eta], \qquad w''(s) = \mathbb{E}[U''(\zeta + s\eta)\eta^2]. \tag{71}$$

PROOF. First, we show that

$$\mathbb{E}[U'(\zeta)\zeta] < \infty. \tag{72}$$

Fix $a > 1$. From (70) and the fact that $U$ is an increasing concave function we deduce that

$$\mathbb{E}[|U(a\zeta)|] < \infty.$$

It follows that

$$\mathbb{E}[U'(a\zeta)\zeta] < \frac{1}{a-1}\mathbb{E}[U(a\zeta) - U(\zeta)] < \infty,$$

which, together with (69), implies (72).

Further, let $0 < b < \frac{1}{K}$. From (13) and (69), we deduce the existence of $c > 0$ such that, for all $|s| \le b$,

$$|U'(\zeta + s\eta)\eta| + |U''(\zeta + s\eta)\eta^2| \le cU'(\zeta)\zeta$$

and, therefore, for all $|s| \le b$ and $|t| \le b$,

$$\left|\frac{U(\zeta + s\eta) - U(\zeta + t\eta)}{s - t}\right| + \left|\frac{U'(\zeta + s\eta) - U'(\zeta + t\eta)}{s - t}\eta\right| \le 2cU'(\zeta)\zeta.$$

The assertion of the lemma now follows from (72) and the Lebesgue theorem on dominated convergence. □

PROOF OF THEOREM 2. As in the statement of the theorem, we fix $x > 0$ and denote $y = u'(x)$. We start with the assertions of item 1. Denote

$$\phi(x) = -\frac{u'(x)}{x}a(x), \tag{73}$$

$$\psi(y) = -\frac{v'(y)}{y}b(y). \tag{74}$$

From (58), we deduce that

$$\phi(x)\psi(y) = -1, \qquad y = u'(x). \tag{75}$$

We claim that

$$u(x+\varepsilon) \geq u(x) + u'(x)\varepsilon + \tfrac{1}{2}\phi(x)\varepsilon^2 + o(\varepsilon^2), \tag{76}$$

$$v(y+\varepsilon) \leq v(y) + v'(y)\varepsilon + \tfrac{1}{2}\psi(y)\varepsilon^2 + o(\varepsilon^2), \tag{77}$$

where we used the standard generic notation $o(\varepsilon)$ for any function $f$ such that $\lim_{\varepsilon\to 0} f(\varepsilon)/|\varepsilon| = 0$.

If $\alpha \in \mathcal{A}^\infty(x)$ then the function

$$w(s) = \mathbb{E}\left[ U\left( g(x)\left( 1 + \frac{s}{x}(1+\alpha) \right) \right) \right],$$

is well defined for sufficiently small $s$ and, by Lemma 4,

$$w'(0) = \frac{1}{x}\mathbb{E}[U'(g(x))g(x)(1+\alpha)] = u'(x)\mathbb{E}_{\mathbb{R}(x)}[1+\alpha] = u'(x), \tag{78}$$

$$w''(0) = \frac{1}{x^2}\mathbb{E}[U''(g(x))g^2(x)(1+\alpha)^2]$$

$$= -\frac{u'(x)}{x}\mathbb{E}_{\mathbb{R}(x)}[A(g(x))(1+\alpha)^2], \tag{79}$$

where $A$ is the relative risk aversion coefficient of $U$ defined in (12) and $\mathbb{R}(x)$ the probability measure introduced in (53). Since $g(x)(1+\frac{s}{x}(1+\alpha)) \in \mathcal{C}(x+s)$ we have

$$u(x+s) \geq w(s).$$

It follows that

$$\liminf_{s\to 0}\left( \frac{u(x+s) - u(x) - u'(x)s}{s^2} \right)$$

$$\geq \liminf_{s\to 0}\left( \frac{w(s) - w(0) - w'(0)s}{s^2} \right)$$

$$= \frac{1}{2}w''(0) = -\frac{u'(x)}{2x}\mathbb{E}_{\mathbb{R}(x)}[A(g(x))(1+\alpha)^2].$$

Taking sup with respect to $\alpha \in \mathcal{A}^\infty(x)$ on the right-hand side of this inequality we deduce that

$$\liminf_{s\to 0}\left( \frac{u(x+s) - u(x) - u'(x)s}{s^2} \right) \geq -\frac{u'(x)}{2x}a(x) = \frac{1}{2}\phi(x),$$

thus proving (76). The proof of (77) is very similar and is omitted here.

Given (76) and (77), the assertions of item 1 are implied by the following result from convex analysis:

LEMMA 5. *Assume that concave functions $u$ and $-v$ are continuously differentiable on $(0, \infty)$ and satisfy the conjugacy relations (9). Let $x > 0$ and $y = u'(x)$. Assume that there are constants $\phi(x)$ and $\psi(y)$ that satisfy (75), (76) and (77).*

*Then, the functions $u$ and $v$ are two-times differentiable at $x$ and $y$, and their respective second derivatives equal $\phi(x)$ and $\psi(y)$.*

PROOF. For sufficiently small $\varepsilon$ we deduce from (9) and (77) that

$$
\begin{aligned}
u(x + \varepsilon) &\leq v(y + \varepsilon\phi(x)) + (x + \varepsilon)(y + \varepsilon\phi(x)) \\
&\leq v(y) + v'(y)\phi(x)\varepsilon + \tfrac{1}{2}\psi(y)(\phi(x)\varepsilon)^2 + (x + \varepsilon)(y + \varepsilon\phi(x)) + o(\varepsilon^2) \\
&= u(x) + u'(x)\varepsilon + \tfrac{1}{2}\phi(x)\varepsilon^2 + o(\varepsilon^2),
\end{aligned}
$$

where at the last step we used (75) and the equality $u(x) = v(y) + xy$.

The last inequality and (76) imply that the function $u$ has the following quadratic expansion at $x > 0$:

$$
u(x + \varepsilon) = u(x) + u'(x)\varepsilon + \tfrac{1}{2}\phi(x)\varepsilon^2 + o(\varepsilon^2).
$$

It is well known (see, e.g., Theorem 5.1.2 in [4]) that for a concave function $u$ the existence of the quadratic expansion at $x$ is equivalent to the two-times differentiability at this point, and that in this case,

$$
u''(x) = \phi(x).
$$

Finally, from (9) and (75), we deduce that $v$ is two-times differentiable at $y = u'(x)$, and

$$
v''(y) = -\frac{1}{u''(x)} = -\frac{1}{\phi(x)} = \psi(y). \qquad \square
$$

The assertions of item 3 are straightforward. Indeed, given the definitions of the random variables $g'(x)$ in (60) and $h'(y)$ in (61), the representations (66) and (67) of the second derivatives of $u$ and $v$ are just reformulations of the equalities in (62) and (63). Further, for $y = u'(x)$, the relation (68) between $g'(x)$ and $h'(y)$ easily follows from the relation (59) between $\alpha(x)$ and $\beta(y)$.

We proceed to the proof of the assertions of item 2. First, we note that convergences (64) and (65) are equivalent. Indeed, given, for example, (64) we deduce from the previous results that

$$
\begin{aligned}
\lim_{\varepsilon \to 0} \frac{h(y + \varepsilon) - h(y)}{\varepsilon} &= \lim_{\varepsilon \to 0} \frac{h(u'(x + \varepsilon)) - h(u'(x))}{u'(x + \varepsilon) - u'(x)} \\
&= \lim_{\varepsilon \to 0} \frac{U'(g(x + \varepsilon)) - U'(g(x))}{\varepsilon} \frac{\varepsilon}{u'(x + \varepsilon) - u'(x)} \\
&= U''(g(x))g'(x)\frac{1}{u''(x)} = h'(y).
\end{aligned}
$$

Now let $(\varepsilon_n)_{n\geq 1}$ be a sequence of real numbers converging to 0. To complete the proof of the theorem it remains to be shown that

$$\lim_{n\to\infty} \frac{g(x+\varepsilon_n)-g(x)}{\varepsilon_n} = g'(x) \tag{80}$$

in probability. By Lemma 3.6 in [7], we have that

$$\lim_{n\to\infty} g(x+\varepsilon_n) = g(x) \tag{81}$$

in probability. In the future it will be convenient for us to assume that, in fact, the convergence in (81) takes place almost surely. Of course, this additional assumption does not restrict any generality.

Define the random variables

$$\zeta = \tfrac{1}{2}\min\left(g(x), \inf_{n\geq 1} g(x+\varepsilon_n)\right),$$

$$\eta = 2\max\left(g(x), \sup_{n\geq 1} g(x+\varepsilon_n)\right),$$

$$\theta = \tfrac{1}{2}\inf_{\zeta\leq t\leq\eta}|U''(t)|.$$

Since the convergence in (81) takes place almost surely, and $g(x) > 0$, we have

$$0 < \zeta < \eta < \infty.$$

From these inequalities and condition (13) on $U$ we deduce that $\theta$ is a strictly positive random variable.

Let $(\alpha_m)_{m\geq 1}$ be a sequence in $\mathcal{A}^\infty(x)$ that converges to the solution $\alpha(x)$ of (56) in $\mathbf{L}_0^2(\mathbb{R}(x))$ and denote

$$g_m = \frac{g(x)}{x}(1+\alpha_m), \qquad m\geq 1.$$

Using the fact that $A(g(x))\leq c_2$, we deduce

$$\lim_{m\to\infty} \mathbb{E}_{\mathbb{R}(x)}[A(g(x))(1+\alpha_m)^2] = \mathbb{E}_{\mathbb{R}(x)}[A(g(x))(1+\alpha(x))^2] = a(x)$$

and, hence, by the definition of $\mathbb{R}(x)$,

$$\lim_{m\to\infty} \mathbb{E}[U''(g(x))g_m^2] = u''(x). \tag{82}$$

For $m\geq 1$, denote by $n_0(m)$ a sufficiently large integer such that, for all $n\geq n_0(m)$,

$$g(x)+\varepsilon_n g_m \in \mathcal{C}(x+\varepsilon_n), \tag{83}$$

$$|\varepsilon_n|(1+|\alpha_m|)\leq \frac{x}{2}. \tag{84}$$

For $n \geq n_0(m)$, we deduce from the Taylor formula and the definition of the random variable $\theta$ that

$$U(g(x) + \varepsilon_n g_m) - U(g(x + \varepsilon_n)) \leq U'(g(x + \varepsilon_n))(g(x) + \varepsilon_n g_m - g(x + \varepsilon_n))$$
$$- \theta(g(x) + \varepsilon_n g_m - g(x + \varepsilon_n))^2.$$

The duality relations (49)–(50) and condition (83) imply that

$$\mathbb{E}[U'(g(x + \varepsilon_n))(g(x) + \varepsilon_n g_m - g(x + \varepsilon_n))] \leq 0.$$

It follows that

$$\mathbb{E}\left[\theta\left(\frac{g(x + \varepsilon_n) - g(x)}{\varepsilon_n} - g_m\right)^2\right] \leq \frac{u(x + \varepsilon_n) - \mathbb{E}[U(g(x) + \varepsilon_n g_m)]}{\varepsilon_n^2}.$$

By Lemma 4 we have

$$\mathbb{E}[U(g(x) + \varepsilon g_m)] = u(x) + u'(x)\varepsilon + \tfrac{1}{2}\mathbb{E}[U''(g(x))g_m^2]\varepsilon^2 + o(\varepsilon^2).$$

Combining this quadratic expansion with convergence (82) we get

$$\lim_{m \to \infty} \lim_{n \to \infty} \frac{u(x + \varepsilon_n) - \mathbb{E}[U(g(x) + \varepsilon_n g_m)]}{\varepsilon_n^2} = 0.$$

Hence,

$$\lim_{m \to \infty} \limsup_{n \to \infty} \mathbb{E}\left[\theta\left(\frac{g(x + \varepsilon_n) - g(x)}{\varepsilon_n} - g_m\right)^2\right] = 0,$$

and the result (80) follows from the facts that $\theta$ is a strictly positive random variable and that the sequence $(g_m)_{m \geq 1}$ converges to $g'(x)$ in probability.

The proof of Theorem 2 is complete. $\square$

**3. Proofs of the main results.** We start with the proof of Lemma 1.

PROOF OF LEMMA 1. The proof is an easy consequence of Lemma 2.

Fix $x > 0$, denote $y = u'(x)$ and let $\mathcal{A}^2(x)$ and $\mathcal{B}^2(y)$ be the sets of the final values of the elements from $\mathcal{M}^2(x)$ and $\mathcal{N}^2(y)$, respectively. That is,

(85)         $\mathcal{A}^2(x) = \{\alpha \in \mathbf{L}^0 : \alpha = M_T \text{ for some } M \in \mathcal{M}^2(x)\},$

(86)         $\mathcal{B}^2(y) = \{\beta \in \mathbf{L}^0 : \beta = N_T \text{ for some } N \in \mathcal{N}^2(y)\}.$

Since $\mathcal{M}^2(x)$ and $\mathcal{N}^2(y)$ are closed linear complementary subspaces in $\mathbf{H}_0^2(\mathbb{R}(x))$, the sets $\mathcal{A}^2(x)$ and $\mathcal{B}^2(y)$ are closed linear complementary subspaces in $\mathbf{L}_0^2(\mathbb{R}(x))$.

Further, we deduce from (13) that

$$c_1 \leq A(X_T(x)) \leq c_2$$

and from (11) and (17) that

$$A(X_T(x))B(Y_T(y)) = 1.$$

Hence, the optimization problems (21) and (22) satisfy the conditions of Lemma 2 and the result follows. $\square$

We proceed now with the proof of the main theorem of the paper.

PROOF OF THEOREM 1. For $x > 0$, denote by $\mathcal{M}^\infty(x)$ the family of semimartingales $M$ such that for some $\varepsilon = \varepsilon(M) > 0$,

$$(87) \qquad X(x)(1 + \varepsilon M) \in \mathcal{X}(x) \quad \text{and} \quad X(x)(1 - \varepsilon M) \in \mathcal{X}(x),$$

where $\mathcal{X}(x)$ is defined in (2). If $M \in \mathcal{M}^\infty(x)$, then $M$ is uniformly bounded and has the initial value 0. Note that the set $\mathcal{M}^\infty(x)$ has the economic interpretation as the set of uniformly bounded wealth processes with initial value 0 under the numéraire $X(x)$.

Similarly, for $y > 0$, denote by $\mathcal{N}^\infty(y)$ the family of semimartingales $N$ such that for some $\varepsilon = \varepsilon(N) > 0$,

$$(88) \qquad Y(y)(1 + \varepsilon N) \in \mathcal{Y}(y) \quad \text{and} \quad Y(y)(1 - \varepsilon N) \in \mathcal{Y}(y),$$

where $\mathcal{Y}(y)$ is defined in (8). If $N \in \mathcal{N}^\infty(y)$, then it is uniformly bounded and $N_0 = 0$.

The following lemma plays the crucial role in the proof.

LEMMA 6. *Assume the conditions of Theorem* 1 *and let* $x > 0$ *and* $y = u'(x)$. *Then, the sets* $\mathcal{M}^\infty(x)$ *and* $\mathcal{N}^\infty(y)$ *belong to* $\mathbf{H}_0^2(\mathbb{R}(x))$, *and their respective closures in* $\mathbf{H}_0^2(\mathbb{R}(x))$ *coincide with the sets* $\mathcal{M}^2(x)$ *and* $\mathcal{N}^2(y)$ *defined in Section* 1.

PROOF. From the definitions of the probability measure $\mathbb{R}(x)$ and the families $\mathcal{M}^\infty(x)$ and $\mathcal{N}^\infty(y)$, we deduce that if $M \in \mathcal{M}^\infty(x)$ and $N \in \mathcal{N}^\infty(y)$, then, for sufficiently small $\varepsilon > 0$, the processes $1 + \varepsilon M$ and $1 - \varepsilon M$, $1 + \varepsilon N$ and $1 - \varepsilon N$, $(1 + \varepsilon M)(1 + \varepsilon N)$ and $(1 + \varepsilon M)(1 - \varepsilon N)$ are supermartingales under $\mathbb{R}(x)$. This clearly implies that the bounded processes $M$, $N$ and $MN$ are martingales under $\mathbb{R}(x)$. Hence, $\mathcal{M}^\infty(x)$ and $\mathcal{N}^\infty(y)$ are orthogonal subspaces in $\mathbf{H}_0^2(\mathbb{R}(x))$.

From (87), we deduce that $\mathcal{M}^\infty(x)$ coincides with the set of *bounded* stochastic integrals with respect to the process $S^{X(x)}$ defined in (19). Hence,

$$\mathcal{M}^\infty(x) \subset \mathcal{M}^2(x).$$

Assumption 2 implies the existence of a strictly positive predictable one-dimensional process $h$ such that the stochastic integral

$$\widetilde{S} = \int h \, dS^{X(x)}$$

is well defined and is a locally bounded process. If $L \in \mathcal{M}^2(x)$, then $L \in \mathbf{H}_0^2(\mathbb{R}(x))$ and there is a predictable process $H$ such that

$$(89) \qquad L = \int H \, d\widetilde{S}.$$

As $\widetilde{S}$ is locally bounded, we deduce that $L$ can be approximated in $\mathbf{H}_0^2(\mathbb{R}(x))$ by bounded stochastic integrals with respect to $\widetilde{S}$ and, hence, by elements from $\mathcal{M}^\infty(x)$. It follows that the closure of $\mathcal{M}^\infty(x)$ in $\mathbf{H}_0^2(\mathbb{R}(x))$ coincides with $\mathcal{M}^2(x)$.

Taking into account orthogonality relations between $\mathcal{M}^\infty(x)$ and $\mathcal{N}^\infty(y)$, we deduce that

$$\mathcal{N}^\infty(y) \subset \mathcal{N}^2(y).$$

To finish the proof, it remains to be shown that any $L \in \mathbf{H}_0^2(\mathbb{R}(x))$ such that $LN$ is a martingale under $\mathbb{R}(x)$ for any $N \in \mathcal{N}^\infty(y)$ is an element of $\mathcal{M}^2(x)$ or, equivalently, has the integral representation (89).

Denote by $\widetilde{\mathcal{Q}}$ the family of equivalent local martingale measures for $\widetilde{S}$ which have bounded densities with respect to $\mathbb{R}(x)$ and by $\widetilde{\mathcal{Z}}$ the family of density processes of these measures. We claim that

$$(90) \qquad Z - 1 \in \mathcal{N}^\infty(y), \qquad Z \in \widetilde{\mathcal{Z}}.$$

Indeed, if $Z \in \widetilde{\mathcal{Z}}$, then for any $X \in \mathcal{X}(1)$ we have that $\frac{X}{X(x)} Z$ is a local martingale under $\mathbb{R}(x)$ and, hence,

$$XZY(y) = \frac{X}{X(x)} Z X(x) Y(y)$$

is a local martingale under $\mathbb{P}$. It follows that

$$(91) \qquad ZY(y) \in \mathcal{Y}(y), \qquad Z \in \widetilde{\mathcal{Z}}.$$

Relation (90) is now implied by (91) and the fact that for any $Z \in \widetilde{\mathcal{Z}}$ there is a sufficiently small $\varepsilon > 0$ such that

$$1 \pm \varepsilon(Z - 1) \in \widetilde{\mathcal{Z}}.$$

From (90) and the assumption that $LN$ is a martingale under $\mathbb{R}(x)$ for all $N \in \mathcal{N}^\infty(y)$, we deduce that $L$ is a martingale under all $\mathbb{Q} \in \widetilde{\mathcal{Q}}$. The integral representation (89) now follows from the well-known result by Jacka; see [5], Theorem 1.1.  □

Given Lemma 6, the proof of the theorem is a rather straightforward consequence of its "abstract version," Theorem 2.

Let $x > 0$ and $y = u'(x)$. By $\mathcal{C}(x)$ and $\mathcal{D}(y)$, we denote the sets of positive random variables which are dominated by the final values of the processes from $\mathcal{X}(x)$ and $\mathcal{Y}(y)$, respectively. That is,

$$\text{(92)} \qquad \mathcal{C}(x) = \{g \in \mathbf{L}^0 : 0 \leq g \leq X_T \text{ for some } X \in \mathcal{X}(x)\},$$

$$\text{(93)} \qquad \mathcal{D}(y) = \{h \in \mathbf{L}^0 : 0 \leq h \leq Y_T \text{ for some } Y \in \mathcal{Y}(y)\}.$$

By $\mathcal{A}^\infty(x)$ and $\mathcal{B}^\infty(y)$, we denote the sets of final values of elements from $\mathcal{M}^\infty(x)$ and $\mathcal{N}^\infty(y)$, respectively. That is,

$$\text{(94)} \qquad \mathcal{A}^\infty(x) = \{\alpha \in \mathbf{L}^\infty : \alpha = M_T \text{ for some } M \in \mathcal{M}^\infty(x)\},$$

$$\text{(95)} \qquad \mathcal{B}^\infty(y) = \{\beta \in \mathbf{L}^\infty : \beta = N_T \text{ for some } N \in \mathcal{N}^\infty(y)\}.$$

With this notation, the value functions $u$ and $v$ defined in (5) and (7) take the form (46) and (47).

According to Proposition 3.1 in [7], the sets $\mathcal{C}(1)$ and $\mathcal{D}(1)$ satisfy conditions (43)–(45). It is easy to see that the sets $\mathcal{C}(x)$ and $\mathcal{D}(y)$ defined in (92) and (93) are related to the sets $\mathcal{A}^\infty(x)$ and $\mathcal{B}^\infty(y)$ defined in (94) and (95) in the same way as the corresponding sets in Section 2, that is, through formulas (51) and (52). Finally, from Lemma 6, we deduce that the respective closures of $\mathcal{A}^\infty(x)$ and $\mathcal{B}^\infty(y)$ in $\mathbf{L}_0^2(\mathbb{R}(x))$ are given by the sets $\mathcal{A}^2(x)$ and $\mathcal{B}^2(y)$ defined in (85) and (86). In particular, these closures are complementary subspaces in $\mathbf{L}_0^2(\mathbb{R}(x))$. Hence, all the assumptions of Theorem 2 are satisfied.

From Theorem 2 we deduce all the assertions of Theorem 1, except the fact that the products $X(x)Y'(y)$, $X'(x)Y(y)$ and $X'(x)Y'(y)$ are martingales under $\mathbb{P}$. However, this result is an immediate consequence of the definitions of the processes $X'(x)$ and $Y'(y)$.   $\square$

## 4. Counterexamples.

This section is devoted to (counter)examples related to our main results. In the first three examples we show that the assertions of Theorem 1 might not hold true if one of the Assumptions 1 or 2 is not satisfied.

EXAMPLE 1. We show that the lower bound

$$c_1 \leq -\frac{xU''(x)}{U'(x)}, \qquad x > 0,$$

in Assumption 1 *cannot* be removed without affecting the existence of the second-order derivative $v''$.

We start by choosing a continuous function $\phi : (0, \infty) \to (0, \infty)$ such that

$$\text{(96)} \qquad \phi(k) = 2^k, \qquad k = 1, 2, \ldots,$$

$$\text{(97)} \qquad \int_0^\infty \phi(s) \, ds = \infty,$$

$$(98) \qquad \int_0^\infty \int_t^\infty \phi(s)\,ds\,dt < \infty$$

and for some $c_2 > 0$,

$$(99) \qquad t\phi(t) > \frac{1}{c_2} \int_t^\infty \phi(s)\,ds, \qquad t > 0.$$

To construct such a function $\phi$, we can start, for example, with the function

$$\psi(t) = \frac{1}{t^{3/2}} e^{-t}, \qquad t > 0,$$

which satisfies (97)–(99), and then modify its values near integers so that (97)–(99) still hold true and, in addition, (96) is satisfied.

We now define

$$V(y) = -\int_0^y \int_t^\infty \phi(s)\,ds\,dt,$$

$$U(x) = \inf_{y>0}\{V(y) + xy\}.$$

Conditions (97) and (98) imply that $U$ satisfies the usual assumptions of a utility function, including Inada conditions (4), and that $U$ and $V$ are bounded:

$$-\infty < V(\infty) = U(0) < V(0) = U(\infty) = 0.$$

In addition, as

$$V'(y) = -\int_y^\infty \phi(t)\,dt, \qquad V''(y) = \phi(y), \qquad y > 0,$$

condition (99) is equivalent to the upper bound in Assumption 1:

$$-\frac{xU''(x)}{U'(x)} < c_2, \qquad x > 0.$$

We now choose a probability space $(\Omega, \mathcal{F}, (\mathcal{F}_t)_{0 \le t \le T}, \mathbb{P})$, where the filtration $(\mathcal{F}_t)_{0 \le t \le T}$ is generated by a Brownian motion $W$. On this probability space, we consider a discrete random variable $\xi$ which takes the values $\frac{1}{2}, 1, 2, 3, \dots, n, \dots$ and such that

$$(100) \qquad \mathbb{E}[\xi] = 1 \quad \text{and} \quad \mathbb{P}[\xi = k] = 2^{-k} \qquad \text{for large } k.$$

It is easy to see that we can construct a continuous stock market (thus trivially satisfying Assumption 2) which is complete and such that the unique martingale measure has the density $\xi$.

For this financial model, the dual value function has the representation

$$v(y) = \mathbb{E}[V(y\xi)], \qquad y > 0.$$

As $V$ is bounded, so is $v$. It follows that $v$ is continuously differentiable (see [8] for a general version of this result) and

$$(101) \qquad v'(y) = \mathbb{E}[\xi V'(y\xi)], \qquad y > 0.$$

Using the Taylor expansion for the function $V$ around $\xi(\omega)$ and (101), we deduce

$$(102) \qquad v(1 + \varepsilon) = v(1) + \varepsilon v'(1) + \tfrac{1}{2}\varepsilon^2 \mathbb{E}[V''(\eta_\varepsilon)\xi^2],$$

where $\eta_\varepsilon(\omega) \in [\xi(\omega), (1 + \varepsilon)\xi(\omega)]$. From (96) and (100), we obtain

$$(103) \qquad \mathbb{E}[V''(\xi)\xi^2] = \mathbb{E}[\phi(\xi)\xi^2] = \infty.$$

Taking into account that $V''(\eta_\varepsilon) \to V''(\xi)$ a.s. we can use (103) and Fatou's lemma to get

$$\lim_{\varepsilon \to 0} \frac{v(1 + \varepsilon) - v(1) - \varepsilon v'(1)}{(1/2)\varepsilon^2} = \infty.$$

Hence $v$ is *not* two-times differentiable at $y = 1$.

EXAMPLE 2. We show that the upper bound

$$-\frac{xU''(x)}{U'(x)} \le c_2, \qquad x > 0,$$

in Assumption 1 is essential for the second-order differentiability of $u$.

In a similar way to Example 1, we can construct a bounded two-times continuously differentiable utility function $U$ such that the lower bound in Assumption 1 holds true, and

$$(104) \qquad -U''(k) = 2^k \qquad \text{for large } k.$$

We consider a one-period financial model with the stock price process $S = (S_0, S_1)$, where $S_0 = 1$ and $S_1$ is a discrete random variable that takes the values $\frac{1}{2}, 1, 2, 3, \ldots, k, \ldots$ and satisfies the following conditions:

$$(105) \qquad \mathbb{P}[S_1 = k] = 2^{-k} \qquad \text{for large } k,$$

$$(106) \qquad \mathbb{E}[U'(S_1)(S_1 - 1)] = 0.$$

Note that, for this model, the set of nonnegative wealth processes with initial value $x$ is given by

$$(107) \qquad \mathcal{X}(x) = \{x + a(S - 1) : a \in [0, 2x]\}.$$

From (106), we deduce that the optimal investment strategy for $x = 1$ is to buy and hold one unit of the stock, that is,

$$(108) \qquad X(1) = S.$$

Indeed, if $X \in \mathcal{X}(1)$, then, by (106) and (107), we have

$$\mathbb{E}[U'(S_1)X_1] = \mathbb{E}[U'(S_1)(1 + a(S_1 - 1))] = \mathbb{E}[U'(S_1)] = \mathbb{E}[U'(S_1)S_1].$$

Using the notation $V$ for the convex conjugate to $U$, we deduce

$$\mathbb{E}[U(X_1)] \le \mathbb{E}[V(U'(S_1)) + U'(S_1)X_1]$$
$$= \mathbb{E}[V(U'(S_1)) + U'(S_1)S_1] = \mathbb{E}[U(S_1)],$$

thus proving (108). We point out that Assumption 2 is satisfied for this model at $x = 1$ because

$$S^{X(1)} = \left( \frac{1}{X(1)}, \frac{S}{X(1)} \right) = \left( \frac{1}{S}, 1 \right)$$

is bounded.

Assumptions (104) and (105) yield

$$(109) \qquad \mathbb{E}[U''(S_1)(1 + a(S_1 - 1))^2] = -\infty, \qquad a \in (-\infty, +\infty).$$

For any $|\varepsilon| < 1$ we know from (107) that

$$u(1 + \varepsilon) = \mathbb{E}[U(1 + \varepsilon + a_\varepsilon(S_1 - 1))]$$
$$= \mathbb{E}\left[ U\left( S_1 + \varepsilon\left( 1 + \frac{a_\varepsilon - 1}{\varepsilon}(S_1 - 1) \right) \right) \right],$$

for some $a_\varepsilon \in [0, 2(1 + \varepsilon)]$. Using the Taylor expansion we obtain

$$u(1 + \varepsilon) = u(1) + \varepsilon u'(1) + \frac{1}{2}\varepsilon^2 \mathbb{E}\left[ U''(\eta_\varepsilon)\left( 1 + \frac{a_\varepsilon - 1}{\varepsilon}(S_1 - 1) \right)^2 \right],$$

where $\eta_\varepsilon$ is a random variable that converges to $S_1$ as $\varepsilon \to 0$. A subsequence argument for $\frac{a_\varepsilon - 1}{\varepsilon}$, together with (109) and Fatou's lemma, now implies that

$$\lim_{\varepsilon \to 0} \frac{u(1 + \varepsilon) - u(1) - \varepsilon u'(1)}{(1/2)\varepsilon^2} = -\infty.$$

Hence, $u$ is *not* two-times differentiable at $x = 1$.

EXAMPLE 3. We show here that if Assumption 2 is violated, then $u$ and $v$ might not be two-times differentiable.

We choose a one-period financial model with one stock, where $S_0 = 1$ and $S_1$ takes the values $2, 1, \frac{1}{2}, \frac{1}{4}, \ldots$ and satisfies

$$(110) \qquad \mathbb{E}\left[ \frac{1}{S_1} \right] = 1 \quad \text{and} \quad \mathbb{E}\left[ \frac{1}{S_1^2} \right] < \infty.$$

We point out that, for this model, the set of nonnegative wealth processes with initial wealth $x$ is given by

$$(111) \qquad \mathcal{X}(x) = \{x + a(S - 1) : a \in [-x, x]\}.$$

We choose two-times continuously differentiable utility function $U$ which is bounded above, which satisfies Assumption 1 and such that

$$(112) \qquad U'(x) = \frac{1}{x}, \qquad x = \frac{1}{2^n}, \qquad n \geq -1,$$

and

$$(113) \qquad \mathbb{E}[U''(S_1)(1 - S_1^2)] = 0.$$

Note that, by Assumption 1 and (112), the random variable $|U''(S_1)S_1^2|$ is bounded and, therefore, the second inequality in (110) implies that

$$\mathbb{E}[|U''(S_1)|(1 + S_1^2)] < \infty.$$

It follows that

$$(114) \qquad \mathbb{E}\left[U''(S_1)\left(\frac{1 + S_1}{2}\right)^2\right] = \sup_a \mathbb{E}[U''(S_1)(1 + a(S_1 - 1))^2] < \infty$$

and, by (113), that the upper bound in (114) is attained at $a = 1/2$.

As in Example 2, we deduce from (110) and (112) that the optimal investment strategy for $x = 1$ is to buy and hold one unit of the stock, that is,

$$X(1) = S.$$

We point out that

$$u'(1) = \mathbb{E}[X_1(1)U'(X_1(1))] = \mathbb{E}\left[S_1 \frac{1}{S_1}\right] = 1,$$

and that the process

$$S^{X(1)} = \left(\frac{1}{X(1)}, \frac{S}{X(1)}\right) = \left(\frac{1}{S}, 1\right)$$

is not sigma-bounded.

For $\varepsilon > 0$, we have $S + \varepsilon\frac{1 + S}{2} \in \mathcal{X}(1 + \varepsilon)$, so

$$(115) \qquad u(1 + \varepsilon) \geq f(\varepsilon) := \mathbb{E}\left[U\left(S_1 + \varepsilon\frac{1 + S_1}{2}\right)\right].$$

From [9], Lemma 1, we deduce that the function $f$ has first and second derivatives from the right at 0 given by

$$(116) \qquad f'_+(0) = \mathbb{E}\left[U'(S_1)\frac{1 + S_1}{2}\right] = 1 = u'(1),$$

$$(117) \qquad f''_+(0) = \mathbb{E}\left[U''(S_1)\left(\frac{1 + S_1}{2}\right)^2\right].$$

Taking into account (114), (115), (116) and (117), we conclude that

$$(118) \qquad \liminf_{\varepsilon \searrow 0} \frac{u(1+\varepsilon) - u(1) - \varepsilon u'(1)}{(1/2)\varepsilon^2} \geq \sup_a \mathbb{E}[U''(S_1)(1 + a(S_1 - 1))^2].$$

From (111), we conclude that, for $\varepsilon > 0$,

$$\frac{X_1(1-\varepsilon) - X_1(1)}{-\varepsilon} = 1 + a_\varepsilon(S_1 - 1) \qquad \text{where } a_\varepsilon \geq 1.$$

Now, according to the Taylor expansion,

$$\frac{u(1-\varepsilon) - u(1) + \varepsilon u'(1)}{(1/2)\varepsilon^2} = \mathbb{E}[U''(\xi_\varepsilon)(1 + a_\varepsilon(S_1 - 1))^2],$$

for $\xi_\varepsilon \in [X_1(1), X_1(1-\varepsilon)]$. Since $U''$ is assumed to be continuous, using a subsequence argument for $a_\varepsilon$ and Fatou's lemma, we obtain

$$(119) \qquad \limsup_{\varepsilon \searrow 0} \frac{u(1-\varepsilon) - u(1) + \varepsilon u'(1)}{(1/2)\varepsilon^2} \leq \sup_{a \geq 1} \mathbb{E}[U''(S_1)(1 + a(S_1 - 1))^2].$$

From (118) and (119), taking into account that the supremum in (114) is strictly attained at $a = \frac{1}{2}$, we conclude that $u$ does not have a second derivative at $x = 1$.

EXAMPLE 4. We construct a very simple financial model that satisfies the conditions of Theorem 1, and such that the derivative processes $X'(x)$ and $Y'(y)$ defined in Section 1 are negative and equal to 0 with a positive probability.

We choose a one-period model such that $S_0 = 1$, and $S_1$ takes the values $\frac{1}{8}, \frac{1}{4}, \frac{1}{2}, 2$ with positive probabilities such that

$$(120) \qquad \mathbb{E}\left[\frac{1}{S_1}\right] = 1.$$

We choose a bounded two-times continuously differentiable utility function $U$ that satisfies Assumption 1 and such that

$$(121) \qquad U'(x) = \frac{1}{x}, \qquad x = \frac{1}{8}, \frac{1}{4}, \frac{1}{2}, 2,$$

$$(122) \qquad \mathbb{E}[U''(S_1)(1 + \tfrac{4}{3}(S_1 - 1))(S_1 - 1)] = 0.$$

Using (120) and (121), we obtain that the optimal investment strategy for the initial capital $x = 1$ is to buy and hold one unit of the stock, that is,

$$X(1) = S.$$

From (122), we deduce that the function

$$f(a) = \mathbb{E}[U''(S_1)(1 + a(S_1 - 1))^2], \qquad a \in (-\infty, +\infty),$$

attains its maximum at $a = 4/3$. It follows from Theorem 1 that the derivative process $X'(1)$ is given by

$$X'(1) = 1 + \tfrac{4}{3}(S - 1).$$

We now have only to see that

$$\mathbb{P}[X_1'(1) = 0] = \mathbb{P}[S_1 = \tfrac{1}{4}] > 0,$$
$$\mathbb{P}[X_1'(1) < 0] = \mathbb{P}[S_1 = \tfrac{1}{8}] > 0.$$

EXAMPLE 5. To motivate the current formulation of Assumption 2 in terms of *sigma-bounded* processes we construct a complete financial model with a *bounded* price process $S$, where, in the case of a logarithmic utility function, the semimartingale $S^{X(x)}$ defined in (19) is not *locally bounded*.

Let $N$ be a Poisson process with intensity 1 defined on a probability space $(\Omega, \mathcal{F}, (\mathcal{F}_t)_{0 \le t \le T}, \mathbb{P})$, where the filtration $(\mathcal{F}_t)_{0 \le t \le T}$ is generated by $N$. We choose a continuous function $\phi \colon (0, \infty) \to (0, \infty)$ such that

$$(123) \qquad\qquad \lim_{t \searrow 0} \phi(t) = \infty,$$

$$(124) \qquad\qquad \int_0^T \phi(t) \, dt = 1/2,$$

and define the process $Z$ as follows:

$$Z_t = 1 + \int_0^t \phi(u) \, d(N_u - u), \qquad 0 \le t \le T.$$

From (124) we conclude that $Z$ is a martingale and that

$$(125) \qquad\qquad Z_t \ge 1/2, \qquad 0 \le t \le T.$$

The process $Z$ is clearly sigma-bounded. However, due to (123) it is not locally bounded. In fact, it is easy to see that the only stopping time $\tau$ such that $Z$ is bounded on $[0, \tau]$ is $\tau = 0$. We define the price process for the stock as

$$S = \frac{1}{Z}.$$

From (125) we deduce that $S$ is a nonnegative price process bounded from above by two. Standard arguments show that $Z$ is the density process of the unique martingale measure and, hence, the model is *complete*.

Consider now the problem of expected utility maximization with the logarithmic utility function $U(x) = \log(x)$. As

$$U'(x) = 1/x, \qquad x > 0,$$

we conclude that the optimal investment strategy for any initial capital is to invest it into the stock, that is,

$$X(x) = xS.$$

The price process of the traded securities (the bond and the stock) under the numéraire $X(x)$ becomes

$$S^{X(x)} = \left( \frac{1}{xS}, \frac{S}{xS} \right) = \frac{1}{x}(Z, 1)$$

and, as we argued above, is sigma-bounded but is *not* locally bounded.

## APPENDIX A: ON SIGMA-BOUNDED SEMIMARTINGALES

In this section we explore the concept of sigma-bounded semimartingales introduced in Definition 1 and present convenient sufficient conditions for the validity of Assumption 2. We start with two simple observations. The first one, for which we skip the proof, is a characterization of locally bounded semimartingales. The second one contains a similar description of sigma-bounded semimartingales. As before, we work on the standard filtered probability space $(\Omega, \mathcal{F}, (\mathcal{F}_t)_{0 \le t \le T}, \mathbb{P})$.

LEMMA 7.   *For a $d$-dimensional semimartingale $R$, the following conditions are equivalent:*

1. *$R$ is locally bounded.*
2. *$R$ has the integral representation*

$$(126) \qquad R_t = R_0 + \int_0^t H_u \, dS_u, \qquad 0 \le t \le T,$$

   *where $S$ is a uniformly bounded semimartingale and $H$ is an increasing, predictable and $S$-integrable process.*
3. *$R$ is dominated by some predictable increasing process $K$, that is,*

$$(127) \qquad \|R_t\| \le K_t, \qquad 0 \le t \le T.$$

LEMMA 8.   *For a $d$-dimensional semimartingale $R$, the following conditions are equivalent:*

1. *$R$ is sigma-bounded.*
2. *$R$ has the integral representation* (126), *where $S$ is a uniformly bounded semimartingale, and $H$ is a predictable and $S$-integrable process.*
3. *$R$ is dominated by some predictable process $K$, that is,* (127) *holds true.*

PROOF. To simplify notation, we assume that all processes appearing in Lemma 8 are one-dimensional. Of course, this assumption does not restrict any generality. The implication $1 \Rightarrow 2$ is straightforward. The implication $2 \Rightarrow 3$ follows from

$$|R_t| \leq |R_{t-}| + |H_t||\Delta S| \leq |R_{t-}| + 2c|H_t| := K_t,$$

where $c$ is a constant dominating the bounded semimartingale $S$. Finally, if (127) holds true for a predictable process $K$, then $|\Delta R| \leq 2K$ and, hence, the stochastic integral $\int \frac{dR}{1+K}$ is locally bounded. This proves $3 \Rightarrow 1$. □

We now state a convenient sufficient condition on the filtered probability space $(\Omega, \mathcal{F}, (\mathcal{F}_t)_{0 \leq t \leq T}, \mathbb{P})$ that implies the sigma-boundedness property for *any* semimartingale. Recall that a local martingale $N$ is called *purely discontinuous* if $NM$ is a local martingale for any continuous local martingale $M$.

ASSUMPTION 4. There is a $d$-dimensional local martingale $M$ such that any bounded, purely discontinuous martingale $N$ is a stochastic integral with respect to $M$, that is,

$$(128) \qquad N_t = N_0 + \int_0^t H_u \, dM_u, \qquad 0 \leq t \leq T,$$

for some predictable and $M$-integrable $H$.

REMARK 3. As the proof of Theorem 3 shows, Assumption 4 is invariant with respect to an equivalent choice of reference probability measure. More precisely, if it holds true under $\mathbb{P}$, then it is also satisfied under any probability measure $\mathbb{Q}$ that is absolutely continuous with respect to $\mathbb{P}$: $\mathbb{Q} \ll \mathbb{P}$. Also, Theorem 3 implies that the integral representation (128) holds true for any (not necessarily bounded) purely discontinuous local martingale $N$.

THEOREM 3. *Assume that Assumption 4 holds true. Then any semimartingale $X$ defined on the filtered probability space $(\Omega, \mathcal{F}, (\mathcal{F}_t)_{0 \leq t \leq T}, \mathbb{P})$ is sigma-bounded.*

The proof of Theorem 3 relies on Proposition 1 below, which is a result of independent interest. For a $d$-dimensional semimartingale $R$ we denote by $\mathcal{I}(R)$ the set of stochastic integrals with respect to $R$, that is,

$$\mathcal{I}(R) = \left\{ X : X = \int H \, dR \text{ for some predictable } H \right\},$$

and by $\mathcal{I}^\infty(R)$ the set of *bounded* elements of $\mathcal{I}(R)$:

$$\mathcal{I}^\infty(R) = \{ X : X \in \mathcal{I}(R) \text{ and is bounded} \}.$$

We use the standard notation $\mathbf{H}^1 = \mathbf{H}^1(\mathbb{P})$ for the space of uniformly integrable martingales under $\mathbb{P}$ such that

$$\|M\|_{\mathbf{H}^1} := \mathbb{E}[M_T^*] < \infty,$$

where $M_T^* = \sup_{0 \leq t \leq T} \|M_t\|$.

PROPOSITION 1.    *Let $R$ be a $d$-dimensional semimartingale. There exists a $d$-dimensional bounded semimartingale $S \in \mathcal{I}^\infty(R)$ such that*

$$\mathcal{I}^\infty(R) = \mathcal{I}^\infty(S).$$

PROOF.    By changing, if necessary, the probability measure $\mathbb{P}$ to an equivalent one, we can assume without any loss of generality that $R$ is a *special* semimartingale, that is,

$$R = R_0 + M + A,$$

where $M$ is a local martingale and $A$ is a predictable process of finite variation. Making the observation that there exists $\varphi > 0$ such that $\overline{M} = \int \varphi \, dM \in \mathbf{H}^1$ and $\overline{A} = \int \varphi \, dA$ has integrable variation, we can also assume that

$$(129) \qquad M \in \mathbf{H}^1 \quad \text{and} \quad \mathbb{E}\left[\int_0^T \|\, dA_t\|\right] < \infty.$$

Let $[M, M] = \sum_{i=1}^d [M^i, M^i]$. We would like to point out that our definition of $[M, M]$ differs from the matrix-valued process often used in the literature. Denoting by $C$ the (one-dimensional) compensator of $[M, M]^{1/2}$, we now define the measure $\mu$ on the predictable $\sigma$-field of $[0, T] \times \Omega$ by

$$d\mu(t, \omega) = (dC_t(\omega) + \|dA_t(\omega)\|) \, d\mathbb{P}(\omega).$$

From (129) we conclude that the measure $\mu$ is finite:

$$\mu([0, T] \times \Omega) < \infty.$$

Denote by $\mathcal{A}$ the set of predictable processes with values in the set of $d \times d$ symmetric and positive semidefinite matrices such that $\int A \, dR$ is well defined and locally bounded. We claim that there is some $\widehat{A} \in \mathcal{A}$ such that

$$(130) \qquad \operatorname{rank}(\widehat{A}) \geq \operatorname{rank}(A) \qquad \forall A \in \mathcal{A},$$

where $\operatorname{rank}(A)$ denotes the *rank* of the matrix $A$ (and the inequality holds $\mu$-a.s.).

Let $(A^n)_{n \geq 1}$ be a sequence in $\mathcal{A}$ such that

$$(131) \qquad \lim_{n \to \infty} \int_{[0,T] \times \Omega} \operatorname{rank}(A^n) \, d\mu = \sup_{A \in \mathcal{A}} \int_{[0,T] \times \Omega} \operatorname{rank}(A) \, d\mu.$$

We claim that the elements of this sequence can be chosen such that

$$\|A_n\| + \|A_n \Delta R\| \leq 1 \tag{132}$$

and

$$\text{rank}(A_{n+1}) \geq \text{rank}(A_n). \tag{133}$$

Indeed, the condition (132) is easy to fulfill, while in order to satisfy (133), it is sufficient to pass from $(A^n)_{n \geq 1}$ to the sequence $(B^n)_{n \geq 1}$ defined by

$$B_1 = A_1,$$

$$B_{n+1} = B_n I_{\{\text{rank}(B_n) \geq \text{rank}(A_{n+1})\}} + A_{n+1} I_{\{\text{rank}(B_n) < \text{rank}(A_{n+1})\}}.$$

Assuming (132) and (133), and denoting by $N$ the predictable process

$$N = \max_{n \geq 1} \text{rank}(A_n)$$

we now define $\widehat{A}$ as

$$\widehat{A} = \sum_{n=1}^{\infty} A_n I_{\{\text{rank}(A_n)=N, \text{rank}(A_k) < N, k=1,2,\ldots,n-1\}}.$$

We deduce that $\widehat{A}$ is an element of $\mathcal{A}$ such that the upper bound on the right-hand side of (131) is attained. As for any $A \in \mathcal{A}$ the sum $A + \widehat{A}$ is in $\mathcal{A}$, and

$$\text{rank}(A + \widehat{A}) \geq \text{rank}(\widehat{A}),$$

it follows that $\widehat{A}$ is also a maximal element in the sense of (130).

We can choose $\widehat{A}$ so that (in addition to (130)) the stochastic integral

$$S = \int \widehat{A} \, dR \tag{134}$$

is bounded (not only locally). Denote by $\mathcal{H}$ the set of $d$-dimensional predictable processes $H$ such that $\int H \, dR$ is locally bounded. To complete the proof, it remains to be shown that any $H \in \mathcal{H}$ admits the representation

$$H = \widehat{A} G \tag{135}$$

for some predictable $d$-dimensional process $G$. Indeed, in this case,

$$\int H \, dR = \int \langle H, dR \rangle = \int \langle \widehat{A} G, dR \rangle = \langle G, \widehat{A} \, dR \rangle = \int G \, dS.$$

To prove (135) we use the fact that any predictable $d$-dimensional process $H$ can be decomposed as

$$H = \widehat{A} G + F, \tag{136}$$

where $G$ and $F$ are predictable $d$-dimensional processes such that $F \in \ker(\widehat{A})$ [$\ker(A)$ denotes the *kernel* of the matrix $A$]. We have to show that for $H \in \mathcal{H}$ the process $F$ in (136) equals zero.

Multiplying, if necessary, both sides of (136) by a strictly positive predictable process, we can assume that $F \in \mathcal{H}$ and $\|F\| \leq 1$. In this case, the matrix $B$ defined by

$$B^{ij} = F^i F^j, \qquad 1 \leq i, j \leq d,$$

belongs to $\mathcal{A}$. Hence, $\widehat{A} + B \in \mathcal{A}$. However, as $F \in \ker(\widehat{A})$

$$\operatorname{rank}(\widehat{A} + B) = \operatorname{rank}(\widehat{A}) + I_{\{F \neq 0\}},$$

and the equality of $F$ to zero follows from the maximality property (130) for $\widehat{A}$.  $\square$

PROOF OF THEOREM 3.  By Proposition 1, we can assume, without loss of generality, that the local martingale $M$ appearing in Assumption 4 is bounded. We can also assume that $M$ is purely discontinuous. We shall maintain these assumptions about $M$ throughout the proof.

Let $Z$ be a bounded martingale orthogonal to $M$, that is, $ZN$ is a martingale and

$$Z = Z^c + Z^d, \qquad Z_0^d = 0,$$

is its decomposition into the sum of the continuous martingale $Z^c$ and the purely discontinuous martingale $Z^d$. As $Z^d$ is locally bounded, it is a stochastic integral with respect to $M$. Since $Z^d$ is orthogonal to $M$, we deduce that $Z^d = 0$.

We have thus proved that any bounded martingale orthogonal to $M$ is continuous. It follows that any purely discontinuous local martingale is a stochastic integral with respect to $M$; see, for example, Theorem 1.1 in [5]. In particular, any purely discontinuous local martingale and, therefore, any special semimartingale, is sigma-bounded.

Since a semimartingale becomes a special semimartingale under an equivalent probability measure, it remains to be shown that Assumption 4 holds true under any $\widetilde{\mathbb{P}} \ll \mathbb{P}$. Let

$$Z_t = \mathbb{E}\left[\frac{d\widetilde{\mathbb{P}}}{d\mathbb{P}} \Big| \mathcal{F}_t\right], \qquad 0 \leq t \leq T,$$

be the density process of $\widetilde{\mathbb{P}}$ with respect to $\mathbb{P}$ and $L = \int dZ/Z_-$ be the stochastic logarithm of $Z$. It is well known that any locally bounded local martingale $\widetilde{N}$ under $\widetilde{\mathbb{P}}$ has the representation

$$\widetilde{N} = N - \langle N, L \rangle, \tag{137}$$

where $N$ is a locally bounded local martingale under $\mathbb{P}$, and $\langle N, L \rangle$ is the predictable covariation between $\mathbb{P}$-local martingales $N$ and $L$. If, in addition, $\widetilde{N}$ is purely discontinuous, then $N$ is also purely discontinuous. Assumption 4 implies the existence of integral representation (128) for $N$ in terms of $M$. It follows that $\widetilde{N}$ has the integral representation (with the same integrand $H$) in terms of $\widetilde{\mathbb{P}}$-local martingale $\widetilde{M}$ defined by

$$\widetilde{M} = M - \langle M, L \rangle.$$

Hence, Assumption 4 is satisfied under $\widetilde{\mathbb{P}}$.   $\square$

The next theorem characterizes local martingales that are sigma-bounded.

THEOREM 4. *For a local martingale $R$, the following assertions are equivalent:*

1. *$R$ is sigma-bounded.*
2. *If $\mathbb{Q}$ is a probability measure such that $\frac{d\mathbb{Q}}{d\mathbb{P}} \in \mathbf{L}^\infty$ and any $X \in \mathcal{I}^\infty(R)$ is a martingale under $\mathbb{Q}$, then any $X \in \mathcal{I}(R) \cap \mathbf{H}^1$ is a martingale under $\mathbb{Q}$.*
3. *The closure of $\mathcal{I}^\infty(R)$ in $\mathbf{H}^1$ coincides with $\mathcal{I}(R) \cap \mathbf{H}^1$.*

PROOF. The implications $1 \Rightarrow 3$ and $3 \Rightarrow 2$ are straightforward. To prove the remaining implication $2 \Rightarrow 1$ we assume that $R \in \mathbf{H}^1$. Of course, this does not restrict any generality.

Let $S$ be the bounded semimartingale (in fact a martingale) given by Proposition 1. Consider a probability measure $\mathbb{Q}$ such that $S$ is a local martingale under $\mathbb{Q}$ and $\frac{d\mathbb{Q}}{d\mathbb{P}} \in \mathbf{L}^\infty$. We have that any element in $\mathcal{I}^\infty(R) = \mathcal{I}^\infty(S)$ is a martingale under $\mathbb{Q}$. According to condition 2, this implies that $R$ is a martingale under $\mathbb{Q}$. From Jacka's theorem (see [5], Theorem 1.1) we deduce that $R$ is a stochastic integral with respect to $S$. The semimartingale $S$ being bounded, we have thus proved that $R$ is sigma-bounded.   $\square$

We conclude this section with an easy corollary of Theorem 3 showing, in particular, that complete financial models satisfy Assumption 2. Hereafter, we consider the financial model with $d$-dimensional price process $S$ and the nonempty family of equivalent martingale measures $\mathcal{Q}$ introduced in Section 1. Recall that the model is called *complete* if any $f \in \mathbf{L}^\infty$ can be represented as the terminal value of a bounded wealth process, that is, $f = X_T$ for some process $X$ such that $X - X_0 \in \mathcal{I}^\infty(S)$. It is well known (see [3]) that the model is complete iff $\mathcal{Q}$ is a singleton.

PROPOSITION 2. *Assume that the financial model is complete. Then any semimartingale $X$ defined on the filtered probability space $(\Omega, \mathcal{F}, (\mathcal{F}_t)_{0 \le t \le T}, \mathbb{P})$ is sigma-bounded.*

Proof.   Denote by $\mathbb{Q}$ the unique element of $\mathcal{Q}$. Assume first that $S$ is a local martingale under $\mathbb{Q}$. As the model is complete, any bounded martingale under $\mathbb{Q}$ belongs to $\mathcal{I}^{\infty}(S)$. Hence Assumption 4 holds true and the result follows from Theorem 3.

In the general case, we use Proposition 1 to find a bounded semimartingale $\widetilde{S} \in \mathcal{I}^{\infty}(S)$ such that

$$\mathcal{I}^{\infty}(S) = \mathcal{I}^{\infty}(\widetilde{S}).$$

We have that $\widetilde{S}$ is a bounded martingale under $\mathbb{Q}$ and that any bounded martingale under $\mathbb{Q}$ belongs to $\mathcal{I}^{\infty}(\widetilde{S})$. The result again follows from Theorem 3.   □

**Acknowledgments.**   We want to thank an Associate Editor and the referees for their comments and remarks. Part of the work was done while Mihai Sîrbu was a graduate student at Carnegie Mellon University.

Department of Mathematical Sciences
Carnegie Mellon University
5000 Forbes Avenue
Pittsburgh, Pennsylvania 15213
USA
E-mail: kramkov@andrew.cmu.edu

Department of Mathematics
Columbia University
2990 Broadway
New York, New York 10027
USA
E-mail: sirbu@math.columbia.edu